\documentclass[12pt]{amsart}
\usepackage[margin=1.in, centering]{geometry}

\usepackage{amssymb,latexsym,amsmath,extarrows, mathrsfs, amsthm}
\usepackage{amsmath}
\usepackage{graphicx}

\usepackage[colorlinks, linkcolor=blue]{hyperref}

\usepackage{graphicx} 

\usepackage{color}
\usepackage{wasysym}
\usepackage{float}

\newtheorem{theorem}{Theorem}[section]
\newtheorem{lemma}[theorem]{Lemma}

\theoremstyle{remark}
\newtheorem*{remark}{\bf  \itshape  \textup{Remark}}

\title[divisor problem in arithmetic progressions]{A new result on the divisor problem in arithmetic progressions modulo a prime power}

\author[Mingxuan Zhong]{Mingxuan Zhong$^{1}$}
\address{1. School of Mathematics and Statistics, Shaanxi Normal University, Xi'an, 710119, Shaanxi, P. R. China}
\email{zhong@snnu.edu.cn}

\author[Tianping Zhang]{Tianping Zhang$^{1,2}$}
\address{2. Research Center for Number Theory and Its Applications, Northwest University, Xi'an, 710127, Shaanxi, P. R. China}
\address{Tianping Zhang is the corresponding author}
\email{tpzhang@snnu.edu.cn}

\allowdisplaybreaks

\date{}

\begin{document}
\def \LLN {\prec \hskip-6pt \prec}

\begin{abstract}
We derive an asymptotic formula for the divisor function $\tau(k)$ in an arithmetic progression $k\equiv a(\bmod \ q)$, uniformly for $q\leq X^{\Delta_{n,l}}$ with $(q,a)=1$. The parameter $\Delta_{n,l}$ is defined as
$$
\Delta_{n,l}=\frac{1-\frac{3}{2^{2^l+2l-3}}}{1-\frac{1}{n2^{l-1}}}.
$$
Specifically, by setting $l=2$, we achieve $\Delta_{n,l}>3/4+5/32$, which surpasses the result obtained by Liu, Shparlinski, and Zhang (2018). Meanwhile, this has also improved upon the result of Wu and Xi (2021). Notably, Hooley, Linnik, and Selberg (1950's) independently established that the asymptotic formula holds for $q\leq X^{2/3-\varepsilon}$. Irving (2015) was the first to surpass the $2/3-$barrier for certain special moduli. We break the classical $3/4-$barrier in the case of prime power moduli and extend the range of $q$. Our main ingredients  borrow from Mangerel's (2021) adaptation of  Mili\'{c}evi\'{c} and Zhang's methodology in dealing with a specific class of weighted Kloosterman sums, rather than adopting Korobov's technique employed by Liu, Shparlinski, and Zhang (2018).

\bigskip

\bigskip


\end{abstract}

\maketitle

\section{Introduction}

\subsection{Background}

Let $q\ge 1$ be a positive integer and $a$ a positive integer such that $(q,a)=1$. For a real number $X\ge 2$, we define
$$
D(q,a;X):=\sum_{n\leq X\atop n\equiv a(q)}\tau(n),
$$
and
$$
E(q,a;X):=D(q,a;X)-\frac{1}{\phi(q)}\sum_{n\leq X\atop (n,q)=1}\tau(n),
$$
where $\phi(q)$ denotes Euler's totient function.

In unpublished works, Hooley, Linnik, and Selberg independently demonstrated that for some constant $\delta>0$ and sufficiently large $X$, the inequality
\begin{align}\label{nontrivial bound}
\max_{(q,a)=1}\left|E(q,a;X)\right|\leq \frac{X^{1-\delta}}{q}
\end{align}
holds uniformly for $q\leq X^{2/3-\varepsilon}$, where $\varepsilon>0$ is an arbitrarily small positive number. This result relies on Weil's bound for Kloosterman sums. It's conjectured that \eqref{nontrivial bound} holds for $q\leq X^{1-\varepsilon}$, but proving this remains challenging, particularly for relatively large values of $q$.

As a compromise, one can consider the mean value of $E(q,a;X)$. For instance, Fouvry \cite{Fouvry1985} proved that for any real number $A>0$ and $a\in \mathbb{Z}$,
$$
\sum_{X^{2/3+\varepsilon}\leq q\leq X^{1-\varepsilon}\atop(q,a)=1}\left|E(q,a;X)\right|\leq X(\log X)^{-A}.
$$
Fouvry and Iwaniec \cite{Fouvry1992} considered the range $X^{2/3-\varepsilon}\leq q\leq X^{2/3+\varepsilon}$ and obtained results that required $q$ to have a good factorization. For the $a-$aspect, Banks, Heath-Brown, and Shparlinski \cite{Banks2005} proved that there exists some constant $\delta>0$ such that
$$
\sum_{1\leq a\leq q\atop (q,a)=1}\left|E(a,q;X)\right|\leq X^{1-\delta}
$$
holds uniformly for $q\leq X^{1-\varepsilon}$. For more related results, we refer to \cite{Blomer2008, Fouvry2015, Friedlander1985 1, Friedlander1985 2, HeathBrown1979, HeathBrown1986, Hooley1957, Lavrik1965, Linnik1961, Wei2016}.

In 2015, Irving \cite{Irving2015} was the first to break $2/3-$barrier and proved that for any real numbers $\varpi, \eta>0$ satisfying $246\varpi+18\eta<1$, there exists a constant $\vartheta>0$, depending on $\varpi$ and $\eta$, such that \eqref{nontrivial bound} holds for any $X^{\eta}-$smooth, square-free $q\leq X^{2/3+\varpi}$ and any integers $a$ with $(q,a)=1$. Khan \cite{Khan2016} considered prime power moduli and obtained the results for $q=p^k, k\ge 7$, where $p$ is a prime. Specifically, there exists some constant $\rho>0$, depending only on $k$, such that \eqref{nontrivial bound} holds for $X^{2/3-\rho}\leq q\leq X^{2/3+\rho}$. Liu, Shparlinski, and Zhang \cite{LiuK2018} also considered prime power moduli and broke the $2/3-$barrier for \eqref{nontrivial bound}, making it uniform in the power of moduli. Wu and Xi \cite{Wu2021} used the theory of exponent pairs and proved that if $q\asymp X^{\theta}$ is square-free and has only prime factors not exceeding $q^{\eta}$ with $\eta>0$ sufficiently small, then \eqref{nontrivial bound} holds for $\theta\leq 2/3+1/232$.

It is natural to ask whether the range of $q$ for which \eqref{nontrivial bound} holds can be extended. Inspired by these conclusions, we investigate the divisor problem modulo a prime power and extend the range of $q$. This result relies on Mili\'{c}evi\'{c} and Zhang's \cite{Milicevic} upper bound for the following type of exponential sums:
$$
\sum_{a\in X}e\left(\frac{f_{T,\boldsymbol{\varepsilon}}(a)}{p^{n}}\right).
$$
We will discuss this issue in Section \ref{sec: proof of exponential sum}.

\subsection{Our results}
We present our results as follows:
\begin{theorem}\label{thm: main theorem}
Let $l\ge 2$, $n>2^{2^{l+1}+6l}$ be fixed integers. Let $q=p^n$ and $a$ be integers, with $p$ an odd prime and $(q,a)=1$. For any sufficiently small constant $\varepsilon>0$, there exists a constant $\Delta_{n,l}$ that depends only on $l$ and $n$, such that
$$
\max_{(q,a)=1}\left|E(q,a;X)\right|\LLN_{p,l}\frac{X^{1-\varepsilon}}{q}
$$
holds uniformly for $q\leq X^{\Delta_{n,l}}$. Specifically, we have
$$
\Delta_{n,l}=\frac{1-\frac{3}{2^{2^l+2l-3}}}{1-\frac{1}{n2^{l-1}}}.
$$
\end{theorem}
Notably, $\Delta_{n,l}$ is a decreasing function of $n$ with a clear lower bound given by
$$\Delta_{n,l}>1-\frac{3}{2^{2^l+2l-3}}.
$$
Furthermore, this lower bound increases with $l$. Even in the worst case where $l=2$, we still have $\Delta_{n,l}>29/32=2/3+23/96=3/4+5/32$, which is an improvement compared to the range of $q$ derived by Liu, Shparlinski and Zhang \cite{LiuK2018}, then later by Wu and Xi \cite{Wu2021}.

For $a\in(\mathbb{Z}/q\mathbb{Z})^{\times}$ and $b\in\mathbb{Z}/q\mathbb{Z}$, we define the classical Kloosterman sums as
\begin{align}\label{classical Kl}
\text{Kl}(a,b;q):=\mathop{{\sum}^*}_{x(q)}e\left(\frac{a\overline{x}+bx}{q}\right).
\end{align}

The proof of Theorem \ref{thm: main theorem} requires the application of a conclusion concerning the following ``weighted Kloosterman sums".

\begin{theorem}\label{thm: weighted Kloosterman sums bound}
Let $M$ be an integer and $\mathcal{A}$ be an interval in $\mathbb{Z}$ with a length not exceeding $K\in\mathbb{R}^{+}$. Suppose that $q=p^n$ with $p>2$ being an odd prime and $n\ge 2$ being a fixed integer. Then, for any fixed positive integer $l\ge 2$ and $n>2^{2^l+3l}$, we have
$$
\sum_{b\in \mathcal{A}}e_{q}(-Mb)\text{Kl}(a,b;q)\LLN_{p,l}Kp^{\frac{n}{2}-\frac{1}{2^{l}}}
$$
uniformly for $K\ge p^{n/2^{2^l +2l-1}+1}$, and any integers $M$ and $a$ satisfying $(q,a)=1$.
\end{theorem}

\begin{remark} From Theorem \ref{thm: main theorem}, it is evident that the range of $q$ increases with $n$. However, to achieve the upper bound of the conjecture, we require $l\rightarrow \infty$, which results in very few $q$ satisfying the condition. Thus, in a sense, the conjecture poses a formidable challenge. \end{remark}

\begin{remark} Compared with previous work of Liu, Shparlinski and Zhang \cite{LiuK2018}, we employ Mili\'{c}evi\'{c} and Zhang's \cite{Milicevic} method to estimate the corresponding exponential sums. Our result establishes a connection between the valid range of the conclusion and the powers of the modulus, but at the cost that as the upper bound’s valid range progressively expands, the number of moduli satisfying the conditions decreases rapidly. \end{remark}

\bigskip

{\bf Notations.} For $x\in\mathbb R$, we always denote $e(x)=e^{2\pi ix}$ and $e_q(x)=e^{\frac{2\pi i x}{q}}$. We use $U\ll V$ or $U=O(V)$ to mean $|U|\leq cV$ for some constant $c>0$. The constant $c$ maybe depend on some parameters $\rho$, in this case, we write $\ll_{\rho}$ or $U=O_{\rho}(V)$. We also use $U\LLN_{\rho}V$, $U\LLN_{\rho,l}V$ to represent $U\leq V\cdot \rho^{o(1)}$ and $U\ll_{l}V\cdot\rho^{o(1)}$ as $\rho\rightarrow\infty$, respectively. If not specified, we do not differentiate between $a\ (q)$ and $a\ (\bmod\ q)$ and they all represent $0\leq a\leq q-1$. For integers $0\leq k< p$, we use $\left(\frac{\cdot}{p^k}\right)$ to denote the Jacobi's symbol. If not specified, we use $\mathop{{\sum}^*_{n}}$ and $\mathop{{\sum}^{\prime}_{n}}$ to represent $\sum_{(n,q)=1}$. We use $\lVert x\rVert$ to represent the distance from $x$ to the nearest integer, and $\lfloor K\rfloor$ ($\left\lceil K\right\rceil$) to represent the lower (upper) rounding function of $K$ respectively.

\bigskip
\section{Kloosterman sums}

In this section, we present some fundamental results concerning the quadratic Kloosterman sums, as defined in Equation \eqref{classical Kl}. Specifically, Lemma \ref{lem: p-adic stationary phase method} corresponds to Equation (12.39) in \cite{Iwaniec2004}, and it is essentially proven using the $p$-adic stationary phase method. Lemma \ref{lem: Kl(a,0,p^n)=0} corresponds to Lemma 3.3 in \cite{LiuK2018}. Lastly, Lemma \ref{lem: Weil's bound on Kl} presents the well-known Weil's bound.

\begin{lemma} \label{lem: p-adic stationary phase method}

Let $q=p^n$, $n\ge 2$ be integers, with $p$ an odd prime. For any $a, b\in\left(\mathbb{Z}/p^n\mathbb{Z}\right)^{\times}$, we show that $\text{Kl}(a,b;p^n)$ vanishes, unless $\left(\frac{ab}{p}\right)=1$ in which case
\begin{align}\label{eq: p-adic stationary phase method}
\text{Kl}(a,b;p^n)=2\left(\frac{l}{p^n}\right)p^{\frac{n}{2}}\text{Re }\varepsilon_{p,n}e\left(\frac{2l}{p^n}\right),
\end{align}
where $l^2\equiv ab(\bmod\ p^n)$, $\varepsilon_{p,n}=1$ if $p^n\equiv 1(\bmod \ 4)$ and $\varepsilon_{p,n}=i$ if $p^n\equiv 3(\bmod \ 4)$.
\end{lemma}

\begin{lemma}
\label{lem: Kl(a,0,p^n)=0}
Let $p$ be an odd prime. For any integer $n\ge 2$, $a\in\left(\mathbb{Z}/p^n\mathbb{Z}\right)^{\times}$ and $b\in\mathbb{Z}/p^n\mathbb{Z}$, if $p\mid b$, then we have $\text{Kl}(a,b;p^n)=0$.
\end{lemma}

\begin{proof}
This is Lemma 3.3 in \cite{LiuK2018}.
\end{proof}

\begin{lemma}\label{lem: Weil's bound on Kl}
Let $q>2$ be an integer, then for any constant $\varepsilon>0$ we have
$$
\mathop{\max}\limits_{\substack{a\in(\mathbb{Z}/q\mathbb{Z})^{\times}\\ b\in\mathbb{Z}/q\mathbb{Z}}}\left|\text{Kl}(a,b;q)\right|\ll_{\varepsilon}q^{1/2+\varepsilon}.
$$
\end{lemma}

\begin{proof}
See \cite{Weil1948}.
\end{proof}

\bigskip

\section{Reduction of the problem}

Let us define the function $E(q,a;X)$ as follows:
$$
E(q,a;X):=\sum_{n\leq X\atop n\equiv a(q)}\tau(n)-\frac{1}{\phi(q)}\sum_{n\leq X\atop (n,q)=1}\tau(n),
$$
with the condition $(q,a)=1$. Starting from the definition of $\tau(n)$, we introduce parameters $U$ and $V$ that satisfy
$$
UV\leq X,\ U\leq X^{1/2},
$$
and $\varepsilon>0$ is sufficiently small. By dyadic subdivision, we focus on a manageable number of sums, specifically $O(\log^2 X)$ sums of the form
$$
E_{1}(q,a;X):=\sum_{u\sim U,v\sim V\atop uv\leq X, uv\equiv a(q)}1-\frac{1}{\phi(q)}\sum_{u\sim U,v\sim V\atop uv\leq X, (uv,q)=1}1,
$$
including an admissible error term.

Next, we further divide the intervals $(U,2U]$, $(V,2V]$ into $O(X^{\delta})$ smaller intervals of lengths $UX^{-\delta}$ and $VX^{-\delta}$, respectively. We denote these intervals as
$$
I_{1}(U_{1}):=\left[U_{1}, U_{1}+UX^{-\delta}\right),\ I_{2}(V_{1}):=\left[V_{1}, V_{1}+VX^{-\delta}\right).
$$
This leads us to consider the case where $U_{1}V_{1}\leq X$, $U_{1}\leq X^{1/2}$. Consequently, $E_{1}(q,a;X)$ is replaced by at most $O(X^{2\delta}\log^2 X)$ sums of the form
$$
E_{2}(q,a;X):=\sum_{u\in I_{1}(U_{1}),v\in I_{2}(V_{1})\atop uv\equiv a(q)}1-\frac{1}{\phi(q)}\sum_{u\in I_{1}(U_{1}),v\in I_{2}(V_{1})\atop (uv,q)=1}1,
$$
with error terms bounded by
$$
\sum_{X<n\leq X+O\left(X^{1-\delta}\right)\atop n\equiv a(q)}\tau(n)\ll\frac{X^{1-\delta+\varepsilon}}{q},\ \frac{1}{\phi(q)}\sum_{X<n\leq X+O\left(X^{1-\delta}\right)\atop (n,q)=1}\tau(n)\ll\frac{X^{1-\delta+\varepsilon}}{q}.
$$

Thus, we obtain the estimate:
\begin{align}\label{eq: Estimation of E(a,q;X)}
E(q,a;X)\ll_{\varepsilon}X^{2\delta+\varepsilon} \max_{U_{1},V_{1}}\left|E_{2}(q,a;X)\right|+\frac{X^{1-\delta+\varepsilon}}{q}.
\end{align}
To proceed with $E_{2}(q,a;X)$, we utilize exponential sums to eliminate the congruence condition and cancel the main term. Specifically:
$$
\sum_{u\in I_{1}(U_{1}),v\in I_{2}(V_{1})\atop uv\equiv a(q)}1=\frac{1}{q}\sum_{k=1}^{q}\left(\sum_{u\in I_{1}(U_{1})\atop (u,q)=1}e_{q}(ak\overline{u})\right)\left(\sum_{v\in I_{2}(V_{1})}e_{q}(-kv)\right).
$$
The term for $k=q$ is
$$
\frac{1}{q}\#\left\{u\in I_{1}(U_{1}), (u,q)=1\right\}\#I_{2}(V_{1}).
$$
Similarly, the second term in $E_{2}(q,a;X)$ can be expressed as
$$
\frac{1}{\phi(q)}\sum_{u\in I_{1}(U_{1}),v\in I_{2}(V_{1})\atop (uv,q)=1}1=\frac{1}{q}\#\left\{u\in I_{1}(U_{1}), (u,q)=1\right\}\#I_{2}(V_{1})+O_{\varepsilon}\left(\frac{X^{1/2}}{q^{1-\varepsilon}}\right)
$$
through a trivial estimation.

Our focus now shifts to bounding
\begin{align}\label{eq: Estimation of E3(a,q;X)}
E_{3}(q,a;X):=&\frac{1}{q}\sum_{k=1}^{q-1}\left|\sum_{u\in I_{1}(U_{1})\atop (u,q)=1}e_{q}(ak\overline{u})\right|\left|\sum_{v\in I_{2}(V_{1})}e_{q}(-kv)\right|\notag\\
=&\frac{1}{q}\sum_{d\mid q}\sum_{k=1\atop (k,q)=d}^{q-1}\left|\sum_{u\in I_{1}(U_{1})\atop (u,q)=1}e_{q}(ak\overline{u})\right|\left|\sum_{v\in I_{2}(V_{1})}e_{q}(-kv)\right|\notag\\
=&\frac{1}{p^n}\sum_{0\leq r<n}\mathop{{\sum}^{\prime}}_{|k|\leq p^{n-r}/2}\left|\sum_{u\in I_{1}(U_{1})\atop \left(u,p^{n-r}\right)=1}e_{p^{n-r}}(ak\overline{u})\right|\left|\sum_{v\in I_{2}(V_{1})}e_{p^{n-r}}(-kv)\right|\notag\\
\ll&\sum_{0\leq r<n}\frac{1}{p^r}\mathop{{\sum}^{\prime}}_{1\leq k\leq p^{n-r}/2}\frac{1}{k}\left|\sum_{u\in I_{1}(U_{1})\atop \left(u,p\right)=1}e_{p^{n-r}}(ak\overline{u})\right|,
\end{align}
provided
$$
\sum_{v\in I_{2}(V_{1})}e_{p^{n-r}}(-kv)\ll \min\left\{VX^{-\delta}, \frac{1}{\lVert k/p^{n-r}\rVert}\right\}.
$$

For the inner sum appearing in \eqref{eq: Estimation of E3(a,q;X)}, we can express it as follows:
\begin{align}\label{eq: inner sum}
&\sum_{u\in I_{1}(U_{1})\atop \left(u,p\right)=1}e_{p^{n-r}}(ak\overline{u})\notag\\
=&\mathop{{\sum}^{\prime}}_{l\left(p^{n-r}\right)}e_{p^{n-r}}\left(ka\overline{l}\right)\sum_{d\in I_{1}(U_{1})\atop d\equiv l\left(p^{n-r}\right)}1\notag\\
=&\frac{1}{p^{n-r}}\sum_{s\left(p^{n-r}\right)}\mathop{{\sum}^{\prime}}_{l\left(p^{n-r}\right)}e_{p^{n-r}}\left(ka\overline{l}+sl\right)\sum_{d\in I_{1}(U_{1})}e_{p^{n-r}}(-sd)\notag\\
=&\frac{1}{p^{n-r}}\sum_{|s|\leq p^{n-r}/2}e_{p^{n-r}}\left(-sU_{1}\right)\text{Kl}\left(ka,s;p^{n-r}\right)g_{p^{n-r}}(s)\notag\\
=&\frac{1}{p^{n-r}}\sum_{1\leq m\leq \frac{p^{n-r}/2}{K}}\sum_{(m-1)K\leq|s|\leq Km}e_{p^{n-r}}\left(-sU_{1}\right)\text{Kl}\left(ka,s;p^{n-r}\right)g_{p^{n-r}}(s),
\end{align}
where we define
$$
g_{w}(s):=e_{w}\left(sU_{1}\right)\sum_{d\in I_{1}\left(U_{1}\right)}e_{w}(-sd),
$$
and set $1\leq K\leq p^{n-r}-1$. Additionally, we observe that
\begin{align}\label{define of gms}
g_{w}(s)=\sum_{d=1}^{U/X^{\delta}}e_{w}(-sd)\ll\min\left\{\frac{U_{1}}{X^{\delta}},\frac{1}{\lVert s/w\rVert}\right\}.
\end{align}

In \eqref{eq: Estimation of E3(a,q;X)}, it is noted that the derived upper bound is acceptable as $p^{r}$ is relatively large, hence our attention can be confined to the case where $p^r$ is small. Introducing a parameter $Z\in\mathbb{Z}$ that will be chosen later, and utilizing \eqref{eq: Estimation of E3(a,q;X)}, Lemma \ref{lem: Kl(a,0,p^n)=0} and Weil's bound for Kloosterman sums, we derive
\begin{align}\label{bound with Z}
&\ \ \sum_{0\leq r<n\atop p^r >Z}\frac{1}{p^r}\mathop{{\sum}^{\prime}}_{1\leq k\leq p^{n-r}/2}\frac{1}{k}\left|\sum_{u\in I_{1}(U_{1})\atop \left(u,p\right)=1}e_{p^{n-r}}(ak\overline{u})\right|\notag\\
&=\sum_{f\mid p^{n}\atop f>Z}\frac{1}{f}\mathop{{\sum}^{\prime}}_{1\leq k\leq p^{n-r}/2}\frac{1}{k}\left|\frac{1}{p^n/f}\mathop{{\sum}^{\prime}}_{|s|\leq p^n/2f}e_{p^n/f}\left(-sU_{1}\right)\text{Kl}\left(ka,s;p^n/f\right)g_{p^n/f}(s)\right|\notag\\
&\LLN_{p}p^{n/2}\cdot Z^{-3/2}.
\end{align}

Next, we define
$$
\kappa:=\max_{1\leq R\leq K}\left|\sum_{s=(K-1)m+1}^{(K-1)m+R}e_{p^{n-r}}\left(-sU_{1}\right)\text{Kl}\left(ka,s;p^{n-r}\right)\right|.
$$
Regarding \eqref{define of gms}, we note that
\begin{align}
\left|g_{w}(s+1)-g_{w}(s)\right|\leq\max_{s\leq t\leq s+1}\left|g_{w}^{\prime}(t)\right|\ll \frac{U_{1}}{wX^{\delta}}\min\left\{\frac{U_{1}}{X^{\delta}},\frac{1}{\lVert s/w\rVert}\right\},\notag
\end{align}
by applying Abel's summation formula to $g^{\prime}_{w}$.

Applying Abel's summation again to \eqref{eq: inner sum}, and combining above with Lemma \ref{lem: Kl(a,0,p^n)=0}, we obtain an upper bound which can be stated as
\begin{align}
&\frac{1}{p^{n-r}}\sum_{1\leq m\leq \frac{p^{n-r}/2}{K}}\kappa\cdot K\max_{(m-1)K\leq|s|\leq Km}\left|g_{p^{n-r}}(s+1)-g_{p^{n-r}}(s)\right|\notag\\
\ll&\frac{U_{1}}{p^{n-r}X^{\delta}}\sum_{1\leq m\leq \frac{p^{n-r}/2}{K}}\frac{\kappa}{m}+K\left(\frac{U_{1}}{p^{n-r}X^\delta}\right)^2\left|\text{Kl}\left(ka,0;p^{n-r}\right)\right|\notag\\
\ll&\frac{U_{1}}{p^{n-r}X^{\delta}}\sum_{1\leq m\leq \frac{p^{n-r}}{K}}\frac{\kappa}{m}.\notag
\end{align}
Finally, back to \eqref{eq: Estimation of E3(a,q;X)}, combining this with \eqref{eq: Estimation of E(a,q;X)} and \eqref{bound with Z}, we arrive at
\begin{align}\label{eq: final es. for E}
E(q,a;X)\LLN_{p} X^{2\delta+\varepsilon} \max_{U_{1},V_{1}}\left|E_{4}(q,a;X)\right|+\frac{q^{1/2}X^{2\delta+\varepsilon}}{Z^{3/2}}+\frac{X^{1-\delta+\varepsilon}}{q},
\end{align}
where
$$
E_{4}(q,a;X):=\frac{U_{1}}{p^n X^{\delta}}\sum_{0\leq r<n\atop p^r\leq Z}\mathop{{\sum}^{\prime}}_{1\leq k\leq p^{n-r}/2}\frac{1}{k}\sum_{1\leq m\leq \frac{p^{n-r}}{K}}\frac{\kappa}{m}.
$$

\bigskip

\section{Proof of Theorem \ref{thm: main theorem}}\label{sec: proof of main thm.}

In this section, we aim to prove Theorem \ref{thm: main theorem} by initially assuming Theorem \ref{thm: weighted Kloosterman sums bound}. To proceed, let $K=X^{\delta}$, Theorem \ref{thm: weighted Kloosterman sums bound} provides us with the bound
$$
\kappa\LLN_{p,l} Kp^{\frac{n-r}{2}-\frac{1}{2^l}},
$$
under the conditions $n-r>2^{2^l+3l}$ and $l\ge 2$. Substituting this into \eqref{eq: final es. for E}, we derive
$$
E_{4}(q,a;X)\LLN_{p,l}\frac{X^{1/2}}{q^{1/2+1/\left(n2^l\right)}}.
$$

Next, we set $\delta=1/2^{2^l+2l-2}$ with $l\ge 2$. Our goal is to ensure that
$$
\frac{X^{1/2}}{q^{1/2+1/\left(n2^l\right)}}\leq \frac{X^{1-3\delta}}{q},
$$
which holds provided
$$
q\leq X^{\Delta_{n,l}},
$$
where
$$
\Delta_{n,l}=\frac{1-\frac{3}{2^{2^l+2l-3}}}{1-\frac{1}{n2^{l-1}}}.
$$
Thus, the first term on the right-hand side of \eqref{eq: final es. for E} is bounded by
$$
\frac{X^{1-3\delta}}{q}X^{2\delta+\varepsilon}\leq \frac{X^{1-\varepsilon}}{q}.
$$

Furthermore, we set $Z=X^{2/3-\varepsilon_{1}-\delta}$ with a sufficiently small $\varepsilon_{1}>0$. Given that $1\leq K<q/Z$ we have $q\ge X^{2/3-\varepsilon_{1}}$. The second term on the right-hand side of \eqref{eq: final es. for E}, which consists of $Z$, is bounded by
$$
\frac{q^{1/2}}{\left(X^{2/3-\varepsilon_{1}-\delta}\right)^{3/2}}X^{2\delta+\varepsilon}\leq \frac{X^{1-\varepsilon}}{q},
$$
provided $q\leq X^{4/3-7\delta/3-4\varepsilon/3-\varepsilon_{1}}$. This condition can be ensured by choosing a large $n$.

Finally, we need to consider the conditions stated in Theorem \ref{thm: weighted Kloosterman sums bound} and the range of $n$. We already have
$$
X^{2/3-\varepsilon_{1}}\leq q\leq X^{\Delta_{n,l}}.
$$
From this, it is straightforward to verify that
$$
X^{\delta}\ge \left(p^{n-r}\right)^{\frac{1}{2^{2^l+2l-1}}+1}.
$$
As for the range of $n$, we only need $$n-r\ge 2^{2^l+3l},$$ with $p^{r}\leq Z$. This can be guaranteed by selecting $n>2^{2^{l+1}+6l}$, which completes the proof of Theorem \ref{thm: main theorem}.

\bigskip

\section{Proof of Theorem \ref{thm: weighted Kloosterman sums bound}}\label{sec: proof of exponential sum}


Let us define the sum $S(N)$ as follows:
$$
S(N):=\sum_{b\in\mathcal{A}(N)}e_q(cb)\text{Kl}(a,b;q),
$$
where $\mathcal{A}(N)=(N,N+K]$ and $|\mathcal{A}(N)|\leq K\in\mathbb{R}^{+}$, $q=p^n$ is the power of an odd prime, and $c$ is an integer. To proceed with our analysis, we will first establish a series of lemmas that will serve as the foundation for the subsequent discussion and proofs.

\subsection{Preliminary lemmas}\label{lemmas for S(N)}
We start with the following lemma. By Weyl's differencing argument, we derive
\begin{lemma}\label{lem: Weyl's difference}
Let $q=p^n$ being the power of an odd prime, $s$, $n\in \mathbb{Z}^{+}$ satisfying $s\leq n$, $H=K/p^s$ and $l\in\mathbb{Z}^{+}$. Then we have
\begin{align}\label{eq: |S(N)|^{2^l}}
\left|S(N)\right|^{2^l}\ll_{l}&q^{2^{l-1}}\sum_{j=1}^{l}\left(p^s\right)^{2^{l-j}}K^{2^l-2^{l-j}}\notag\\
&+K^{2^l-l-1}\left(p^s\right)^l \sum_{0< \left|h_{1}\right|,\dots,\left|h_{l}\right|\leq H}\left|\sum_{k\in\mathcal{A}\left(N;h_{1},\dots, h_{l}\right)}\prod_{J\subseteq \left\{1,\dots ,l\right\}}\text{Kl}\left(a,k+\sum_{j\in J}p^sh_{j};q\right)\right|,
\end{align}
where $\mathcal{A}\left(N;h_{1},\dots, h_{l}\right):=(N-p^sh_{1}-\dots-p^sh_{l},N+K-p^sh_{1}-\dots-p^sh_{l}]$.
\end{lemma}

\begin{proof} Actually, this is a conclusion analogous to Lemma 4.3 in \cite{Irving2015}, but focused specifically on the modulo $p^n$ setting.

For an integer $J^{\prime}\geq 1$ and some integers $s_{j}\in \mathbb{Z}$, $j=1,2,\dots,J^{\prime}$, we define a more general sum as
$$
S(J^{\prime},N):=\sum_{k\in \mathcal{A}(N)}e_{q}(ck)\prod_{j=1}^{J^{\prime}}\text{Kl}\left(a,k+s_{j};q\right).
$$
We then introduce the sequence $a_{k}$ defined by
$$
a_{k}:=
\begin{cases}
e_{q}(kc)\prod_{j=1}^{J^{\prime}}\text{Kl}\left(a,k+s_{j};q\right), & k\in \mathcal{A}(N)\\
0,& k\notin \mathcal{A}(N)
\end{cases}
$$
to transform $S(J^{\prime},N)$ into a bilinear form:
\begin{align}
S(J^{\prime},N)=&\sum_{k}a_{k}\notag\\
=&\frac{1}{H}\sum_{h=1}^{H}\sum_{k}a_{k+p^s h}\notag\\
=&\frac{1}{H}\sum_{k}\sum_{h=1}^{H}a_{k+p^s h},\notag
\end{align}
where $H=K/p^s$ and $1\leq s\leq n$. If $k+p^{s}h\in \mathcal{A}(N)$, then
\begin{align}
a_{k+p^s h}=&e_{q}\left(\left(k+p^s h\right)c\right)\prod_{j=1}^{J^{\prime}}\text{Kl}\left(a,k+p^s h+s_{j};q\right)\notag\\
=&e_{q}(kc)e_{q}\left(p^s hc\right)\prod_{j=1}^{J^{\prime}}\text{Kl}\left(a,k+p^s h+s_{j};q\right).\notag
\end{align}


Applying Cauchy-Schwarz inequality, we find that
\begin{align}
&H^2 \left|S(J^{\prime},N)\right|^2 \notag\\
\leq&K\sum_{k}\left|\sum_{h=1\atop k+p^sh\in\mathcal{A}(N)}^{H}e_{q}\left(p^shc\right)\prod_{j=1}^{J^{\prime}}\text{Kl}\left(a,k+p^s h+s_{j};q\right)\right|^2\notag\\
\leq& K\sum_{h_{1},h_{2}=1}^{H}\left|\sum_{k\atop k+p^sh_{1},k+p^sh_{2}\in\mathcal{A}(N)}\prod_{j=1}^{J^{\prime}}\text{Kl}\left(a,k+p^s h_{1}+s_{j};q\right)\text{Kl}\left(a,k+p^s h_{2}+s_{j};q\right)\right|\notag\\
=&K\sum_{h_{1},h_{2}=1}^{H}\left|\sum_{k\in\mathcal{A}(N)\atop k+p^s(h_{1}-h_{2})\in\mathcal{A}(N)}\prod_{j=1}^{J^{\prime}}\text{Kl}\left(a,k+s_{j};q\right)\text{Kl}\left(a,k+p^s (h_{2}-h_{1})+s_{j};q\right)\right|\notag\\
\leq&KH\sum_{|h|\leq H}\left|\sum_{k\in\mathcal{A}(N)\atop k+p^sh\in\mathcal{A}(N)}\prod_{j=1}^{J^{\prime}}\text{Kl}\left(a,k+s_{j};q\right)\text{Kl}\left(a,k+p^sh+s_{j};q\right)\right|,\notag
\end{align}
which further implies
\begin{align}\label{eq: induction for l=1}
&\left|S(J^{\prime},N)\right|^2\notag\\
\leq&p^s\left(Kq^{J^{\prime}}+\sum_{0<|h|\leq H}\left|\sum_{k\in\mathcal{A}(N;h)}\prod_{j=1}^{J^{\prime}}\text{Kl}\left(a,k+s_{j};q\right)\text{Kl}\left(a,k+p^sh+s_{j};q\right)\right|\right),
\end{align}
where $\mathcal{A}(N;h):=\mathcal{A}(N)-p^sh=(N-p^sh,N+K-p^sh]$, the same as below.

The rest of the proof can be obtained by induction. When $l=1$, the result follows from \eqref{eq: induction for l=1}. Assuming that
\begin{align}
\left|S(N)\right|^{2^{l-1}}\ll_{l}&q^{2^{l-2}}\sum_{j=1}^{l-1}\left(p^s\right)^{2^{l-1-j}}K^{2^{l-1}-2^{l-1-j}}+K^{2^{l-1}-l}\left(p^s\right)^{l-1} \notag\\
&\times\sum_{0< \left|h_{1}\right|,\dots,\left|h_{l-1}\right|\leq H}\left|\sum_{k\in\mathcal{A}\left(N;h_{1},\dots, h_{l-1}\right)}\prod_{J\subseteq\left\{1,\dots, l-1\right\}}\text{Kl}\left(a,k+\sum_{j\in J}p^sh_{j};q\right)\right|,\notag
\end{align}
squaring both sides yields
\begin{align}
\left|S(N)\right|^{2^l}\ll_{l}&q^{2^{l-1}}\sum_{j=1}^{l-1}\left(p^s\right)^{2^{l-j}}K^{2^l-2^{l-j}}+K^{2^l-2l}\left(p^s\right)^{2l-2}\notag\\
&\times\left( \sum_{0< \left|h_{1}\right|,\dots,\left|h_{l-1}\right|\leq H}\left|\sum_{k\in\mathcal{A}\left(N;h_{1},\dots, h_{l-1}\right)}\prod_{J\subseteq\left\{1,\dots, l-1\right\}}\text{Kl}\left(a,k+\sum_{j\in J}p^sh_{j};q\right)\right|\right)^2.\notag
\end{align}
Using Cauchy-Schwarz inequality for the inner sums one by one, we obtain
\begin{align}
&\ \ \ \sum_{0< \left|h_{1}\right|,\dots,\left|h_{l-1}\right|\leq H}\left|\sum_{k\in\mathcal{A}\left(N;h_{1},\dots, h_{l-1}\right)}\prod_{J\subseteq\left\{1,\dots, l-1\right\}}\text{Kl}\left(a,k+\sum_{j\in J}p^sh_{j};q\right)\right|\notag\\
&\leq\sum_{0< \left|h_{1}\right|,\dots,\left|h_{l-2}\right|\leq H} \left(\sum_{0<\left|h_{l-1}\right|\leq H}1^2\right)^{\frac{1}{2}}\notag\\
&\ \ \ \times\left(\sum_{0<\left|h_{l-1}\right|\leq H}\left|\sum_{k\in\mathcal{A}\left(N;h_{1},\dots, h_{l-1}\right)}\prod_{J\subseteq\left\{1,\dots, l-1\right\}}\text{Kl}\left(a,k+\sum_{j\in J}p^sh_{j};q\right)\right|^2\right)^{\frac{1}{2}}\notag\\
&\leq H^{\frac{1}{2}}\sum_{0< \left|h_{1}\right|,\dots,\left|h_{l-3}\right|\leq H} \left(\sum_{0<\left|h_{l-2}\right|\leq H}1^2\right)^{\frac{1}{2}}\notag\\
&\ \ \ \times\left(\sum_{0<\left|h_{l-2}\right|, \left|h_{l-1}\right|\leq H}\left|\sum_{k\in\mathcal{A}\left(N;h_{1},\dots, h_{l-1}\right)}\prod_{J\subseteq\left\{1,\dots, l-1\right\}}\text{Kl}\left(a,k+\sum_{j\in J}p^sh_{j};q\right)\right|^2\right)^{\frac{1}{2}}\notag\\
&\leq\cdots\notag\\
&\cdots\notag\\
&\leq H^{\frac{l-1}{2}}\left(\sum_{0< \left|h_{1}\right|,\dots,\left|h_{l-1}\right|\leq H}\left|\sum_{k\in\mathcal{A}\left(N;h_{1},\dots, h_{l-1}\right)}\prod_{J\subseteq\left\{1,\dots, l-1\right\}}\text{Kl}\left(a,k+\sum_{j\in J}p^sh_{j};q\right)\right|^2\right)^{\frac{1}{2}},\notag
\end{align}
which leads to
\begin{align}\label{eq: induction for l=l}
\left|S(N)\right|^{2^l}\ll_{l}&q^{2^{l-1}}\sum_{j=1}^{l-1}\left(p^s\right)^{2^{l-j}}K^{2^l-2^{l-j}}+K^{2^l-l-1}\left(p^s\right)^{l-1}\notag\\
&\times\sum_{0< \left|h_{1}\right|,\dots,\left|h_{l-1}\right|\leq H}\left|\sum_{k\in\mathcal{A}\left(N;h_{1},\dots, h_{l-1}\right)}\prod_{J\subseteq\left\{1,\dots, l-1\right\}}\text{Kl}\left(a,k+\sum_{j\in J}p^sh_{j};q\right)\right|^2.
\end{align}

In \eqref{eq: induction for l=1}, taking $J^{\prime}=2^{l-1}$, $c=0$ and $s_{j}=\sum_{0\leq i\leq j-1}p^sh_{i}$ then we have
\begin{align}
&\left|\sum_{k\in\mathcal{A}\left(N;h_{1},\dots ,h_{l-1}\right)}\prod_{J\subseteq\left\{1,\dots, l-1\right\}}\text{Kl}\left(a,k+\sum_{j\in J}p^sh_{j};q\right)\right|^2\notag\\
\leq&p^sKq^{2^{l-1}}+p^s\sum_{0<|h_{l-1}|\leq H}\left|\sum_{k\in\mathcal{A}\left(N;h_{1},\dots, h_{l}\right)}\prod_{J\subseteq\left\{1,\dots, l-1\right\}}\text{Kl}\left(a,k+\sum_{j\in J}p^sh_{j};q\right)\right.\notag\\
&\times \left.\text{Kl}\left(a,k+\sum_{j\in J}p^sh_{j}+p^sh_{l};q\right)\right|\notag\\
=&p^sKq^{2^{l-1}}+p^s\sum_{0<|h_{l}|\leq H}\left|\sum_{k\in\mathcal{A}\left(N;h_{1},\dots, h_{l}\right)}\prod_{J\subseteq\left\{1,\dots, l\right\}}\text{Kl}\left(a,k+\sum_{j\in J}p^sh_{j};q\right)\right|.\notag
\end{align}
Now Lemma \ref{lem: Weyl's difference} follows from substituting the above into \eqref{eq: induction for l=l}.\end{proof}

The next lemma is an estimation of exponential sums with argument function $f_{T,\boldsymbol{\varepsilon}}$ and set $X$ (see \cite{Milicevic} for the definition), proved by Mili\'{c}evi\'{c} and Zhang.
\begin{lemma}\label{lem: Estimates on complete exponential sums}
For every $M>0$, there exists constants $\delta_{i}=\delta_{i}(M)>0, (i=1,2,3)$ such that for every $n\in \mathbb{N}$, $T\subseteq \mathbb{Z}/p^{n}\mathbb{Z}$, every set $X\subseteq (\mathbb{Z}/p^{n}\mathbb{Z})^{[T]}$ that is invariant under translations by $p^{\lfloor \delta_{2}n\rfloor}\mathbb{Z}/p^{n}\mathbb{Z}$, $\boldsymbol{\varepsilon}=(\varepsilon_{\tau})_{\tau\in T} \in\mathbb{Z}^{T}$ with $\lVert \varepsilon\rVert_{1}=\sum_{\tau\in T}|\varepsilon_{\tau}|\leq M$, either
$$
\sum_{a\in X}e\left(\frac{f_{T,\boldsymbol{\varepsilon}}(a)}{p^{n}}\right)\ll_{|T|}p^{(1-\delta_{1})n}
$$
or
$$
\tau\equiv \tau^{\prime}(\bmod \ p^{\lfloor \delta_{3}(\lfloor\delta_{2} n\rfloor-\rho)\rfloor})
$$
for some $\tau, \tau^{\prime}\in T$ with $\tau\neq \tau^{\prime}$ and $\varepsilon_{\tau}, \varepsilon_{\tau^{\prime}}\not\equiv 0(\bmod \ p^{n})$.
\end{lemma}

\begin{proof} See Theorem 4, Proposition 8 and 9 in \cite{Milicevic}.\end{proof}

For trivial cases in Lemma \ref{lem: Estimates on complete exponential sums}, we introduce the following results. Lemma \ref{lem: trivial case=0} is the second case of dealing with Lemma \ref{lem: Estimates on complete exponential sums}, while Lemmas \ref{mixed character sum} and \ref{cyclotomic polynomial} are used to deal with the case of $\boldsymbol{\varepsilon}=\boldsymbol{0}$.

Let $L$ be a positive integer, define $$T:=\left\{H_{I}=\sum_{i\in I}h_{i}:\ I\subseteq\left\{1,\dots,L\right\}\right\},\ \boldsymbol{h}:=\left\{h_{1},h_{2},\dots,h_{L}\right\}\in \mathbb{Z}^{L},$$ and $$\mu(\tau):=\left|\left\{I\subseteq\left\{1,\dots,L\right\}:\  H_{I}\equiv \tau(\bmod\ p^{v})\right\}\right|.$$ Then we have

\begin{lemma}\label{lem: trivial case=0}
Let $c\in\left(0,2^{-2^L}\right]$ be a real number and let $Q_{j}\in \mathbb{Z}$, $j=0,1,\dots,L$. Assume that $d\mid Q_{0}$, $K/Q_{j}\ge Q_{0}^{2c}$ for all $1\leq j\leq L$. Then the number of tuples $\boldsymbol{h}\in \mathbb{Z}^{L}$ with $1\leq |h_{j}|\leq K/Q_{j}$ for all $j$, such that for each $p^{v}\Vert d$ and
$$
\min_{I,J\subseteq \left\{1,\dots, L\right\}\atop H_{I}\neq H_{J}}\left|H_{I}- H_{J}\right|_{p}<p^{-cv}
$$
is $\ll 2^{L}\tau(d)^{2L-1}d^{-c}K^{L}/(Q_{1}\cdots Q_{L})$.
\end{lemma}

\begin{proof} See Lemma 5.3 in \cite{Mangerel2021}.\end{proof}

\begin{lemma}\label{mixed character sum}
Let $p$ be an odd prime and $C\in \mathbb{Z}/p\mathbb{Z}$. Let $A\in (\mathbb{Z}/p\mathbb{Z})^{\times}$ and $T\subseteq \mathbb{Z}/p\mathbb{Z}$, then we have
$$
\sum_{d(p)\atop A(d+\tau)\in (\mathbb{Z}/p\mathbb{Z})^{\times 2}, \forall \tau\in T}e_{p}(dC)=2^{-|T|}p\textbf{1}_{C\equiv 0(p)}+O_{|T|}\left(p^{1/2}\right).
$$
\end{lemma}

\begin{proof} Let $\chi_{2}$ be Legendre's symbol. Then we obtain
\begin{align}
&\sum_{d(p)\atop A(d+\tau)\in (\mathbb{Z}/p\mathbb{Z})^{\times 2}, \forall \tau\in T}e_{p}(dC)\notag\\
=&\sum_{d(p)}e_{p}(dC)\prod_{\tau\in T}\textbf{1}_{A(d+\tau)\in (\mathbb{Z}/p\mathbb{Z})^{\times 2}}\notag\\
=&2^{-|T|}\sum_{\boldsymbol{j}\in\left\{0,1\right\}^{|T|}}\chi_{2}(A)^{t(\boldsymbol{j})}\sum_{d(p)}\chi_{2}\left(\prod_{\tau\in T}(d+\tau)^{j_{\tau}}\right)e_{p}(dC),\notag
\end{align}
where $t(\boldsymbol{j})=\sum_{\tau\in T}j_{\tau}$. The main term comes from $\boldsymbol{j}=\boldsymbol{0}$, which we state as
$$
2^{-|T|}\sum_{d(p)}e_{p}(dC)+O_{|T|}(1)=2^{-|T|}p\textbf{1}_{C\equiv 0(p)}+O_{|T|}(1).
$$

When $\boldsymbol{j}\neq \boldsymbol{0}$, from Weil's bound in \cite{Weil1948} and Corollary 11.24 in \cite{Iwaniec2004} we have
\begin{align}
&2^{-|T|}\sum_{\boldsymbol{j}\in\left\{0,1\right\}^{|T|}\atop \boldsymbol{j}\neq \boldsymbol{0}}\chi_{2}(A)^{t(\boldsymbol{j})}\sum_{d(p)}\chi_{2}\left(\prod_{\tau\in T}(d+\tau)^{j_{\tau}}\right)e_{p}(dC)\notag\\
\ll&_{|T|}\max_{\boldsymbol{j}\neq \boldsymbol{0}}\left|\sum_{d(p)}\chi_{2}\left(\prod_{\tau\in T}(d+\tau)^{j_{\tau}}\right)e_{p}(dC)\right|\ll p^{1/2}.\notag
\end{align}
\end{proof}

\begin{lemma}\label{cyclotomic polynomial}
Let $p$ be an odd prime, $p\mid Q_{0}\in \mathbb{Z}$ and let $\boldsymbol{h}\in \mathbb{Z}^{L}$. If $2\mid \mu(\tau)$ for all $\tau\in T$, then there exists some $1\leq i\leq L$ such that $p\mid h_{i}$.
\end{lemma}

\begin{proof} Assume $p\nmid h_{i},\ \forall i$. We begin our analysis with the summation
$$
\sum_{I\subseteq\left\{1,\dots,L\right\}}e_{p}(aH_{I}),\ p\nmid a.
$$
First, we express this summation in a product form. Specifically, we have
$$
\sum_{I\subseteq\left\{1,\dots,L\right\}}e_{p}(aH_{I})=\prod_{1\leq i\leq L}\left(1+e_{p}(ah_{i})\right).
$$
Next, we consider the product over $a\in(\mathbb{Z}/p\mathbb{Z})^{\times}$. By doing so, we obtain
$$
\prod_{a\in(\mathbb{Z}/p\mathbb{Z})^{\times}}\prod_{1\leq i\leq L}\left(1+e_{p}(ah_{i})\right)=\prod_{b\in(\mathbb{Z}/p\mathbb{Z})^{\times}}\left(1+e_{p}(b)\right)^{L}=1,
$$
where we have applied the property of the cyclotomic polynomial of order $p$.

On the other hand, we can also represent the original summation in a different form. Specifically, we have
$$
\sum_{I\subseteq\left\{1,\dots,L\right\}}e_{p}(aH_{I})=\sum_{\tau\in T}\mu(\tau)e_{p}(a\tau).
$$
Now, let us consider the case where $2\mid \mu(\tau)$. In this case, we can find some $m_{\tau}\in \mathbb{Z}$ such that
$$
\sum_{\tau\in T}\mu(\tau)e_{p}(a\tau)=2\sum_{\tau\in T}m_{\tau}e_{p}(a\tau).
$$
Similarly, by taking the product over $a\in(\mathbb{Z}/p\mathbb{Z})^{\times}$, we get
$$
2^{L-1}A_{p}:=2^{L-1}\prod_{a\in(\mathbb{Z}/p\mathbb{Z})^{\times}}\sum_{\tau\in T}m_{\tau}e_{p}(a\tau).
$$
In view of Lemma 5.1 in \cite{Mangerel2021}, we know that $A_{p}$ is a rational integer. However, this leads to a contradiction, as it implies that $2\mid 1$. Consequently, our initial assumption must be incorrect, and the lemma is proved.\end{proof}

\subsection{Upper bound for $S(N)$}\label{estimation for S(N)}
To use the results presented by Mili{\'c}evi{\'c} and Zhang \cite{Milicevic} (i.e. Lemma \ref{lem: Estimates on complete exponential sums}), it is necessary to reformulate $S(N)$ in a manner that aligns with the structure outlined in Lemma \ref{lem: Estimates on complete exponential sums}. Consequently, we will follow a similar approach as outlined in Section 4 of \cite{Mangerel2021} to establish the required exponential sums.

We begin with the inner sum from Lemma \ref{lem: Weyl's difference}, and define
\begin{align}
S(l,N,\boldsymbol{h}):=\sum_{k\in\mathcal{A}\left(N;h_{1},\dots, h_{l}\right)}\prod_{J\subseteq\left\{1,\dots, l\right\}}\text{Kl}\left(a,k+\sum_{j\in J}p^sh_{j};q\right).\notag
\end{align}

Next, we introduce the following notations
$$
\mu(\tau)=\mu_{\boldsymbol{h}}(\tau):=\left|\left\{J\subseteq\left\{1,\dots, l\right\}:\ \sum_{j\in J}p^sh_{j}\equiv \tau(\bmod\  q)\right\}\right|,
$$
and
$$
T=T_{\mu}:=\left\{\tau\in\mathbb{Z}/p^n\mathbb{Z}:\mu(\tau)\ge 1\right\}.
$$
Using these notations, we can rewrite $S(l,N,\boldsymbol{h})$ as
\begin{align}\label{eq: define of S(l,N)}
S(l,N,\boldsymbol{h})=\sum_{k\in\mathcal{A}\left(N;h_{1},\dots, h_{l}\right)}\prod_{\tau\in T}\text{Kl}\left(a,k+\tau;q\right)^{\mu(\tau)}.
\end{align}

To proceed, let $s(r)$ be a solution to the congruence $s(r)^2\equiv r(\bmod\ p^{\lfloor n/2\rfloor})$ and $\mu_{0}$ be a primitive square root modulo $p^{\lfloor n/2\rfloor}$. Based on Hensel's lemma, if $x^2\equiv r\left(\bmod\ p^{\lfloor n/2\rfloor}\right)$, then $x\equiv s(r)\mu_{0}^j\ \left(\bmod\ p^{\lfloor n/2\rfloor}\right)$ (this is a lift of a primitive square root modulo $p$) for some $0\leq j\leq d_{p}-1$, where $d_{p}=2$.

Now, we apply Lemma \ref{lem: p-adic stationary phase method} to transform the Kloosterman sums into exponential sums. As outlined in the lemma, the summation consists of two terms, denoted as $y_{1}$ and $y_{2}$. From the definitions, we have
$$
y_{1}\equiv s(r)(\bmod\ p), \ y_{2}\equiv s(r)\mu_{0}(\bmod\ p).
$$
It follows that $y_{1}\equiv -y_{2} (\bmod\  p)$, which implies
$$
\left(\frac{y_{1}}{p}\right)=\left(\frac{-1}{p}\right)\left(\frac{y_{2}}{p}\right).
$$
Furthermore, when $n$ is odd, the Kloosterman sums can be expressed as
$$
\text{Kl}\left(a,k+\tau;p^n\right)=
\begin{cases}
\pm\left(\frac{a}{p}\right)p^{\frac{n}{2}}\sum_{0\leq j\leq d_{p}-1}e_{p^n}\left(2as(\overline{a}(k+\tau))\mu_{0}^j\right),\ p\equiv 1(\bmod\ 4),\\
\mp\left(\frac{a}{p}\right) i p^{\frac{n}{2}}\sum_{0\leq j\leq d_{p}-1}(-1)^{j}e_{p^n}\left(2as(\overline{a}(k+\tau))\mu_{0}^j\right),\ p\equiv 3(\bmod\ 4).\\
\end{cases}
$$
When $n$ is even, Kloosterman sums can be simply expressed as
$$
\text{Kl}\left(a,k+\tau;p^n\right)=p^{\frac{n}{2}}\sum_{0\leq j\leq d_{p}-1}e_{p^n}\left(2as(\overline{a}(k+\tau))\mu_{0}^j\right).
$$
It is worth noting that the signs in the above expressions do not affect the final result in $S(N)$. In fact, for odd $n$ with $p\equiv 1(\bmod\ 4)$ or $p\equiv 3(\bmod\ 4)$, and even $n$ the upper bounds remain the same.

For instance, when $n$ is odd and $p\equiv 3(\bmod\ 4)$, substituting the above into \eqref{eq: define of S(l,N)} yields
\begin{align}
S(l,N,\boldsymbol{h})=&\left(\frac{a}{p}\right)^{2^l}i^{2^l}\left(p^{\frac{n}{2}}\right)^{2^l}\sum_{\boldsymbol{j}=(\boldsymbol{j}(\tau))_{\tau}\in U}(-1)^{\sum_{\tau\in T}\sum_{1\leq j\leq \mu(\tau)}j_{i}(\tau)}\notag\\
& \times\sum_{k\in\mathcal{A}\left(N;h_{1},\dots, h_{l}\right)}e_{p^n}\left(2a\sum_{\tau\in T}s(\overline{a}(k+\tau))\left(\sum_{1\leq i\leq \mu(\tau)}\mu_{0}^{j_{i}(\tau)}\right)\right),\notag
\end{align}
where
$$
U(\tau):=\left\{0,\dots, d_{p}-1\right\}^{\mu(\tau)},\ U=\prod_{\tau\in T}U(\tau).
$$

To simplify further, we introduce
$$
f_{T,\boldsymbol{\varepsilon}}(k):=2a\sum_{\tau\in T}\boldsymbol{\varepsilon}_{\tau}s(\overline{a}(k+\tau)),
$$
$$
\phi(\boldsymbol{\varepsilon}):=\left|\left\{\boldsymbol{j}\in U:\boldsymbol{\varepsilon}_{\tau}=\sum_{1\leq i\leq \mu(\tau)}\mu_{0}^{j_{i}(\tau)} \right\}\right|,
$$
for each $\boldsymbol{\varepsilon}\in (\mathbb{Z}/p^{n}\mathbb{Z})^{|T|}$. With these definitions, \eqref{eq: define of S(l,N)} can be reduced to
$$
S(l,N,\boldsymbol{h})\leq\left(p^{\frac{n}{2}}\right)^{2^l}\sum_{\boldsymbol{\varepsilon}\in\left(\mathbb{Z}/p^n\mathbb{Z}\right)^{|T|}}
\phi(\boldsymbol{\varepsilon})\left|\sum_{k\in\mathcal{A}\left(N;h_{1},\dots, h_{l}\right)}e_{p^n}\left(f_{T,\boldsymbol{\varepsilon}}(k)\right)\right|.
$$

It is important to observe that, for odd $n$ with $p\equiv 1(\bmod\ 4)$ or even $n$, $S(l,N,\boldsymbol{h})$ remains the same upper bound.

It is beneficial for our proof process to complete the inner sum as follows
\begin{align}\label{eq: inner sum in S(l,N)}
\sum_{k\in\mathcal{A}\left(N;h_{1},\dots, h_{l}\right)}e_{p^n}\left(f_{T,\boldsymbol{\varepsilon}}(k)\right)=&\frac{1}{p^n}\sum_{a(p^n)}\sum_{k(p^n)}e_{p^n}\left(f_{T,\boldsymbol{\varepsilon}}(k)+ak\right)
\sum_{b\in\mathcal{A}\left(N;h_{1},\dots, h_{l}\right)}e_{p^n}(-ab)\notag\\
=&\frac{1}{p^n}\sum_{1\leq |a|\leq p^n/2}\sum_{k(p^n)}e_{p^n}\left(f_{T,\boldsymbol{\varepsilon}}(k)+ak\right)
\sum_{b\in\mathcal{A}\left(N;h_{1},\dots, h_{l}\right)}e_{p^n}(-ab)\notag\\
&+\frac{\lfloor K\rfloor}{p^n}\sum_{k\left(p^n\right)}e_{p^n}\left(f_{T,\boldsymbol{\varepsilon}}(k)\right)\notag\\
=&\frac{1}{p^n}\sum_{1\leq m\leq p^n/2K^{\prime}}\sum_{K^{\prime}(m-1)<|r|\leq K^{\prime}m}\sum_{k(p^n)}e_{p^n}\left(f_{T,\boldsymbol{\varepsilon}}(k)+rk\right)\notag\\
&\times\sum_{b\in\mathcal{A}\left(N;h_{1}\dots h_{l}\right)}e_{p^n}(-rb)+\frac{\lfloor K\rfloor}{p^n}\sum_{k\left(p^n\right)}e_{p^n}\left(f_{T,\boldsymbol{\varepsilon}}(k)\right)\notag\\
=:&\frac{1}{p^n}\sum_{1\leq m\leq p^n/2K^{\prime}}S+\frac{\lfloor K\rfloor}{p^n}\sum_{k\left(p^n\right)}e_{p^n}\left(f_{T,\boldsymbol{\varepsilon}}(k)\right),
\end{align}
where $1\leq 2K^{\prime}\leq p^n$.

\begin{remark} Here, we once again establish a bilinear sum for the summation over $a$. This step is because there may exist trivial terms in the sum, which we wish to minimize their impact in the estimation. This can be ensured by taking an appropriately small value of $K^{\prime}$. \end{remark}

Next, let's define
$$
g(r):=\sum_{b\in\mathcal{A}\left(N;h_{1},\dots, h_{l}\right)}e_{p^n}(-rb),
$$
and apply partial summation to $S$. We find that
\begin{align}
S\ll&\sum_{r=K^{\prime}(m-1)+1}^{K^{\prime}m}\left|\sum_{k(p^n)}e_{p^n}\left(f_{T,\boldsymbol{\varepsilon}}(k)+rk\right)\right|\cdot K^{\prime}\max_{K^{\prime}(m-1)<r\leq K^{'}m}\left|g(r+1)-g(r)\right|\notag\\
\ll&\sum_{r=K^{\prime}(m-1)+1}^{K^{\prime}m}\left|\sum_{k(p^n)}e_{p^n}\left(f_{T,\boldsymbol{\varepsilon}}(k)+rk\right)\right|\cdot K^{\prime}\max_{K^{\prime}(m-1)<r\leq K^{\prime}m}\max_{r\leq t\leq r+1}\left|g^{\prime}(t)\right|\notag\\
\ll&\sum_{r=K^{\prime}(m-1)+1}^{K^{\prime}m}\left|\sum_{k(p^n)}e_{p^n}\left(f_{T,\boldsymbol{\varepsilon}}(k)+rk\right)\right|\cdot K^{\prime}\frac{K}{p^n}\min\left\{K,\frac{1}{ K^{\prime}m/p^n}\right\}\notag\\
\ll&\sum_{r=K^{\prime}(m-1)+1}^{K^{\prime}m}\left|\sum_{k(p^n)}e_{p^n}\left(f_{T,\boldsymbol{\varepsilon}}(k)+rk\right)\right|\cdot \frac{K}{m}.\notag
\end{align}
Now, by choosing $K^{\prime}=q^{\varepsilon}$ and substituting it into \eqref{eq: inner sum in S(l,N)}, we obtain
\begin{align}\label{eq: complete the inner sum in S(l,N)}
\sum_{k\in\mathcal{A}\left(N;h_{1},\dots ,h_{l}\right)}e_{p^n}\left(f_{T,\boldsymbol{\varepsilon}}(k)\right)\LLN_{p,l}&\frac{K}{p^n}\left(\sum_{1\leq m\leq p^n/2q^{\varepsilon}}\frac{1}{m}\sum_{r=q^{\varepsilon}(m-1)+1}^{q^{\varepsilon}m}\left|\sum_{k(p^n)}e_{p^n}\left(f_{T,\boldsymbol{\varepsilon}}(k)+rk\right) \right|\right.\notag\\
&+\left.\left|\sum_{k(p^n)}e_{p^n}\left(f_{T,\boldsymbol{\varepsilon}}(k)\right)\right|\right).
\end{align}
Furthermore, this leads to
\begin{align}\label{eq: upper bound of S(l,N)}
S(l,N,\boldsymbol{h})\LLN_{p,l}&\left(p^{\frac{n}{2}}\right)^{2^l}\frac{K}{p^n}\sum_{1\leq m\leq p^n/2q^{\varepsilon}}\frac{1}{m}\sum_{r=q^{\varepsilon}(m-1)+1}^{q^{\varepsilon}m}\sum_{\boldsymbol{\varepsilon}\in\left(\mathbb{Z}/p^n\mathbb{Z}\right)^{|T|}}
\phi(\boldsymbol{\varepsilon})\left|\sum_{k(p^n)}e_{p^n}\left(f_{T,\boldsymbol{\varepsilon}}(k)+rk\right) \right|\notag\\
&+\left(p^{\frac{n}{2}}\right)^{2^l}\frac{K}{p^n}\sum_{\boldsymbol{\varepsilon}\in\left(\mathbb{Z}/p^n\mathbb{Z}\right)^{|T|}}
\phi(\boldsymbol{\varepsilon})\left|\sum_{k(p^n)}e_{p^n}\left(f_{T,\boldsymbol{\varepsilon}}(k)\right)\right|.
\end{align}

We can now finally use Lemma \ref{lem: Estimates on complete exponential sums} to analyze the inner sums in \eqref{eq: upper bound of S(l,N)}. Specifically, by setting
$$
X:=\left\{k(\bmod \ p^n):\ k\equiv d(\bmod\  p)\right\},
$$
we observe that this set is $p^{\lfloor \delta_{2}n\rfloor}\mathbb{Z}/p^{n}\mathbb{Z}-$invariant. Consequently, indicated in Lemma \ref{lem: Estimates on complete exponential sums}, if $\boldsymbol{\varepsilon}\neq \boldsymbol{0}$ and $n>2^{2^l+3l}$, then the following estimate holds
\begin{align}
\sum_{k(p^n)\atop k\equiv d(p)}e_{p^n}\left(f_{T,\boldsymbol{\varepsilon}}(k)+rk\right)\ll p^{(1-\delta_{1})n}.\notag
\end{align}
Alternatively, there exist at least two distinct $\tau,\tau^{\prime}\in T$, such that $\varepsilon_{\tau},\varepsilon_{\tau^{\prime}}\not\equiv 0(\bmod\ p^n)$ and
$$
\tau\equiv\tau^{\prime}\left(\bmod\ p^{\lfloor\delta_{3}\left(\lfloor\delta_{2}n\rfloor-\rho\right)\rfloor}\right),
$$
where $\delta_{i}>0,i=1,2,3$, $\rho>0$ and $d(\bmod\ p)$ satisfies $\overline{a}(d+\tau)\in\left(\left(\mathbb{Z}/p^n\mathbb{Z}\right)^{\times}\right)^{2}$.

More precisely, we select $\delta_{1}=\delta_{2}=2^{-|T|}$, $\delta_{3}=\tbinom{|T|}{2}^{-1}$ and
$$
\rho=\textbf{1}_{p\leq 2|T|-1}+\left\lceil\frac{\log\left(2\cdot 2^{l}\right)}{\log p}\right\rceil.
$$
The validity of these values is ensured by Proposition 4.8 in \cite{Mangerel2021}. Furthermore, the value of $\rho$ comprises two parts, which can be derived through similar calculations as those in Lemma 4.4 and Lemma 4.9 in \cite{Mangerel2021}, respectively.

Given that $|T|\leq2^l$, we can simplify the estimate to
\begin{align}\label{eq: estimation of the inner sum of S(l,N)}
\sum_{k(p^n)\atop k\equiv d(p)}e_{p^n}\left(f_{T,\boldsymbol{\varepsilon}}(k)+rk\right)\ll p^{(1-2^{-|T|})n},
\end{align}
and express another case as
\begin{align}\label{eq: condition of Milicevic and Zhang's estimation}
\tau\equiv\tau^{\prime}\left(\bmod\ p^{D(l)}\right),
\end{align}
with
$$D(l)=\left\lfloor\tbinom{|T|}{2}^{-1}\left(\left\lfloor2^{-|T|}n\right\rfloor-\textbf{1}_{p\leq 2\times|T|-1}-\left\lceil\frac{\log\left(2\cdot 2^{l}\right)}{\log p}\right\rceil\right)\right\rfloor.$$
Now, using \eqref{eq: estimation of the inner sum of S(l,N)}, we can deduce
\begin{align}\label{eq: upper bound of S(l,N) 2}
S(l,N,\boldsymbol{h})\LLN_{q,l}Kp^{n2^{l-1}-n/2^{2^l}+1},
\end{align}
with $n>2^{2^l+3l}$ and $\boldsymbol{\varepsilon}\neq \boldsymbol{0}$.

We now need to address some trivial cases, namely $\boldsymbol{\varepsilon}= \boldsymbol{0}$ and \eqref{eq: condition of Milicevic and Zhang's estimation} holds.

Firstly, recall the notation of $\tau$, if $\boldsymbol{h}$ fulfills the condition \eqref{eq: condition of Milicevic and Zhang's estimation}, as stated in Lemma \ref{lem: trivial case=0}, since $p^n \lVert d$ we can set $K=K/p^{s}$, $d=Q_{0}=p^n$, $Q_{i}=1,1\leq i\leq l$ and $c=1/2^{2^l+2l}<1/2^{2l}$. This implies that the number of solutions $\boldsymbol{h}$ with $1\leq |h_{i}|\leq K/p^{s}$ satisfying
$$
\min_{I,J\subseteq \left\{1,\dots, l\right\}\atop H_{I}\neq H_{J}}\left|H_{I}-H_{J}\right|_{p}<p^{-cn}
$$
is less than
$$
2^l\tau(p^n)^{2l-1}p^{-cn}\left(K/p^s\right)^l,
$$
where $H_{I}=\sum_{i\in I}h_{i}$. This conclusion holds by choosing $K\ge p^{n/2^{2^l +2l-1}+1}$, and $l\ge2$ later. Substituting above and \eqref{eq: upper bound of S(l,N) 2} into \eqref{eq: |S(N)|^{2^l}}, we obtain the upper bound for the inner sum in Lemma \ref{lem: Weyl's difference} as
\begin{align}\label{eq: inner sum of |S(N)|^{2^l}}
\sum_{0< \left|h_{1}\right|,\dots,\left|h_{l}\right|\leq H}S(l,N,\boldsymbol{h})\LLN_{p,l}KH^{l}p^{n2^{l-1}-\frac{n}{2^{2^l+2l}}}.
\end{align}

Secondly, let's consider the case where $\boldsymbol{\varepsilon}= \boldsymbol{0}$. Recall the definition of $\phi(\boldsymbol{\varepsilon})$. If $\phi(\boldsymbol{\varepsilon})\neq 0$, there exists tuples $\boldsymbol{j}\in U$ such that
$$
0=\alpha_{\tau}+\beta_{\tau}\mu_{0}
$$
for some nonnegative integers $\alpha_{\tau}$ and $\beta_{\tau}$. Additionally, we have
$$
\alpha_{\tau}+\beta_{\tau}=\mu(\tau)
$$
for all $\tau\in T$. Setting $r:=\left\lceil\frac{\log\left(2\cdot 2^{l}\right)}{\log p}\right\rceil\leq n/2$, it is evident that
$$
p^{r}\mid \alpha_{\tau}+\beta_{\tau}\mu_{0}=\alpha_{\tau}-\beta_{\tau}+\beta_{\tau}(1+\mu_{0})\Rightarrow p^r\mid(\alpha_{\tau}-\beta_{\tau}),
$$
and $|\alpha_{\tau}-\beta_{\tau}|\leq 2\cdot2^{l}\leq p^r$. Consequently, we deduce that $\alpha_{\tau}=\beta_{\tau}$ and $2\mid \mu(\tau)$.

By incorporating Lemma \ref{mixed character sum}, we know that in this case
\begin{align}
\phi(\boldsymbol{0})\left|\sum_{k(p^n)}e_{p^n}\left(f_{T,\boldsymbol{0}}(k)+rk\right)\right|=0,\notag
\end{align}
unless $p^{n-1}\mid r$ and $2\mid \mu(\tau)$ for all $\tau\in T$. Furthermore,
\begin{align}
\phi(\boldsymbol{0})\left|\sum_{k(p^n)}e_{p^n}\left(f_{T,\boldsymbol{\varepsilon}}(k)+rk\right)\right|\ll p^{n-1}\cdot
p^{\frac{1}{2}}, \text{ if } r\neq 0 \text{ and } 2\mid \mu(\tau).\notag
\end{align}
Using Lemma \ref{cyclotomic polynomial}, the condition $2\mid \mu(\tau)$ can be translated into a constraint on $\boldsymbol{h}$ stating that
$$
 \exists 1\leq i\leq l\text{ such that } p\mid h_{i}\Leftrightarrow p\mid\prod_{1\leq i\leq l}h_{i}.
$$
This leads us to
\begin{align}\label{eq: upper bound of S(l,N) 3}
S(l,N,\boldsymbol{h})&\LLN_{p,l}Kp^{n2^{l-1}}\textbf{1}_{p\mid \prod_{i}h_{i}}\left(1+\sum_{1\leq m\leq p^n/2q^{\varepsilon}}\frac{1}{m}\sum_{r=q^{\varepsilon}(m-1)+1}^{q^{\varepsilon}m}p^{-\frac{1}{2}}\cdot \textbf{1}_{p^{n-1}\mid r}\right)\notag\\
&\ll Kp^{n2^{l-1}}\textbf{1}_{p\mid \prod_{i}h_{i}}.
\end{align}
Substituting \eqref{eq: upper bound of S(l,N) 3} into \eqref{eq: |S(N)|^{2^l}} again, we obtain that in this case, the upper bound for the inner sum in Lemma \ref{lem: Weyl's difference} as
$$
\sum_{0< \left|h_{1}\right|,\dots,\left|h_{l}\right|\leq H}S(l,N,\boldsymbol{h})\LLN_{p,l} KH^{l}p^{n2^{l-1}-1}.
$$
Combining this with \eqref{eq: inner sum of |S(N)|^{2^l}}, we conclude that
\begin{align}
\sum_{0< \left|h_{1}\right|,\dots,\left|h_{l}\right|\leq H}S(l,N,\boldsymbol{h})\LLN_{p,l}&\max\left\{KH^{l}p^{n2^{l-1}-1},KH^{l}p^{n2^{l-1}-\frac{n}{2^{2^l+2l}}}\right\}\notag\\
\LLN_{p,l}&KH^{l}p^{n2^{l-1}-1}\notag
\end{align}
for $n>2^{2^l+3l}$.

Thus taking $s=1$ in \eqref{eq: |S(N)|^{2^l}} yields
\begin{align}\label{eq: upper bound for S(N)}
S(N)\LLN_{p,l}Kp^{\frac{n}{2}}\sum_{j=1}^{l}\left(\frac{p}{K}\right)^{\frac{1}{2^j}}+Kp^{\frac{n}{2}-\frac{1}{2^l}}\LLN Kp^{\frac{n}{2}-\frac{1}{2^l}},
\end{align}
where $n>2^{2^l+3l}$ and $K\ge p^{n/2^{2^l +2l-1}+1}$. This ends the proof of Theorem \ref{thm: weighted Kloosterman sums bound}.

\bigskip

\section*{Acknowledgements}

This work is supported by the National Natural Science Foundation of China (Nos. 12271320, 11871317), and the Natural Science Basic Research Plan for Distinguished Young Scholars in Shaanxi Province of China (No. 2021JC-29). The authors sincerely thank Professor I. E. Shparlinski (UNSW) and Professor Kui Liu (QDU) for their assistance in this project. The authors also sincerely appreciate the valuable suggestions from the anonymous reviewers, which have effectively enhanced the quality of this manuscript.

\baselineskip=0.9\normalbaselineskip


{\small
\begin{thebibliography}{99}

\bibitem{Banks2005}Banks W. D., Heath-Brown D. R. and Shparlinski I. E., On the average value of divisor sums in arithmetic progression, IMRN, 2005, 1: 1--25.

\bibitem{Blomer2008}Blomer V., The average value of divisor sums in arithmetic progressions, Q. J. Math., 2008, 59: 275--286.

\bibitem{Fouvry1985}Fouvry \'{E}, Sur le probl\ {e}me des diviseurs de Titchmarsh, J. Reine Angew. Math., 1985, 357: 51--76.

\bibitem{Fouvry1992}Fouvry \'{E}, Iwaniec H., The divisor function over arithmetic progressions, Acta Arith., 1992, 61(3): 271--87.

\bibitem{Fouvry2015}Fouvry \'{E}, Kowalaski E. and Michel P., On the exponent of distribution of the ternary divisor function, Mathematika, 2015, 61(1): 121--144.

\bibitem{Friedlander1985 1}Friedlander J. B., Iwaniec H., The divisor problem for arithmetic progression, Acta Arith., 1985, 45: 273--277.

\bibitem{Friedlander1985 2}Friedlander J. B., Iwaniec H., Incomplete Kloosterman sums and a divisor problem, Ann. Math., 1985, 121(2): 319--355.

\bibitem{HeathBrown1979}Heath-Brown D. R., The fourth power moment of the Riemann zeta function, Proc. Lond. Math. Soc., 1979, 38(3): 385--422.

\bibitem{HeathBrown1986}Heath-Brown D. R., The divisor function $d_{3}(n)$ in arithmetic progression, Acta Arith., 1986, 47: 29--56.

\bibitem{Hooley1957}Hooley C., An asymptotic formula in the theory of numbers, Proc. Lond. Math. Soc., 1957, 7(3): 396--413.

\bibitem{Irving2015}Irving A., The divisor function in arithmetic progression to smooth moduli, IMRN, 2015, 15: 6675--6698.

\bibitem{Iwaniec2004}Iwaniec H., Kowalski E., Analytic Number Theory. Amer. Math. Soc. Colloq. Publ., AMS, Providence, RI, 2004, 53, xii+615 pp.

\bibitem{Khan2016}Khan R., The divisor function in arithmetic progressions modulo prime powers, Mathematika, 2016, 62: 898--908.

\bibitem{Lavrik1965}Larvik A. F., On the problem of divisors in segments of arithmetical progression(in Russian), Dokl. Akad. Nauk SSSR, 1965, 164(6): 1232--1234.

\bibitem{Linnik1961}Linnik Yu. V., All large numbers are sums of a prime and two squares (a problem of Hardy and Littlewood) \uppercase\expandafter{\romannumeral2}, Mat. Sb., 1961, 53(95): 3--38.

\bibitem{LiuK2018}Liu K., Shparlinski I. E. and Zhang T. P., Divisor problem in arithmetic progressions modulo a prime power, Adv. Math., 2018, 325: 459--481.

\bibitem{Mangerel2021}Mangerel A. P., Squarefree integers in arithmetic progressions to smooth moduli, Forum. Math. Sigma., 2021, 9(Paper No. e72): 47 pp.

\bibitem{Milicevic}Mili{\'c}evi{\'c} D., Zhang S., Distribution of Kloosterman paths to high prime power moduli, Trans. Amer. Math. Soc. Ser. B, 2023, 10: 636--669.

\bibitem{Wei2016}Wei F., Xue B. Q. and Zhang Y. T., General divisor functions in arithmetic progressions to large moduli. Sci. China Math., 2016, 59(9), 1663--1668.

\bibitem{Weil1948}Weil A., On some exponential sums, Proc. Nat. Acad. Sci. U.S.A, 1948, 34: 204--207.

\bibitem{Wu2021}Wu J., Xi, P., Arithmetic exponent pairs for algebraic trace functions and applications, Algebra Number Theory, 2021, 15(9): 2123--2172.

\end{thebibliography}
}
\end{document}